\documentstyle[12pt]{article}
\begin{document}
\newtheorem{Def}{Definition}[section]
\newtheorem{thm}{Theorem}[section]
\newtheorem{lem}{Lemma}[section]
\newtheorem{rem}{Remark}[section]
\newtheorem{prop}{Proposition}[section]
\newtheorem{cor}{Corollary}[section]
\def\av{{\int \hspace{-2.25ex}-} }
\title
{Some multi-valued solutions to Monge-Amp\`ere equations}
\author{ L. Caffarelli\thanks{Partially
 supported by
 NFS grant DMS-0140388 and G-37-X71-G4}\\
Department of Mathematics
\\
The University of Texas
\\
Austin, TX 78712
\\ 
\\ YanYan Li
\thanks{Partially
 supported by
    NSF grant DMS-0401118.  
    }
    \\ Department of Mathematics
\\
  Rutgers University\\
     110 Frelinghuysen Rd. \\
        Piscataway, NJ 08854}
\date{}
\maketitle
\input { amssym.def}

\bigskip

\setcounter{section}{-1}
\section{Introduction}
In this paper we construct several types of
multi-valued solutions to the  Monge-Amp\`ere equation
in higher dimensions.  Recently, there has been considerable
interest in understanding the behavior of a metric generated by
a solution of the  Monge-Amp\`ere equation at a singularity.  See for instance
Loftin, Yau and Zaslow
\cite{LYZ} and 
Leung \cite{Le}.   To explain our results let us refer first to the theory of 
multi-valued harmonic functions.
Typical two dimensional examples
of multi-valued harmonic functions are
$$
f_1(z)=Re(z^{\frac 1k})
$$
and
$$
f_2(z)=arg(z).
$$
The first one, $f_1$, is finitely valued,
since it repeats itself once you
have gone $k-$times around the origin.  In the case  
of $f_2$, each time we go around the origin, the function increases by
$2\pi$.
Multi-valued harmonic functions have been studied by
G.V. Evans (\cite{E1}, \cite{E2} and \cite{E3}),
H. Lewy (\cite{L}) and L. Caffarelli (\cite{C4} and \cite{C5}), which
have inspired the present paper.

\section{Finitely valued solutions of the
Monge-Amp\`ere equation}
The geometric situation is the following:
Let $D\subset \Bbb R^n$, $n\ge 2$,
be  a bounded  strictly convex open set
with smooth boundary $\partial D$,
and let
$\Sigma\subset D$ be homeomorphic in $\Bbb R^n$ to an 
$(n-1)-$dimensional closed disc, i.e., there exists a
homeomorphism $\psi: \Bbb R^n\to \Bbb R^n$ such that
$\psi(\Sigma)$ is  an
$(n-1)-$dimensional closed disc.  Let $\Gamma=\partial
\Sigma$, the boundary of $\Sigma$.  Thus $\Gamma$ is homeomorphic to 
an $(n-2)-$dimensional sphere for $n\ge 3$.
In $\Bbb R^3$, $\Gamma$ is a curve ``spanned'' by a  disc.

Let 
$$
M=(D \setminus \Gamma)\times \Bbb Z,
$$
denote a covering of $D\setminus \Gamma$
with the following standard parameterization: Fixing an $x^*\in
D\setminus \Sigma$, and connecting $x^*$ by a smooth curve 
in $D\setminus \Gamma$ to a point $x$ in $D\setminus \Gamma$.
If the curve goes through $\Sigma$ $m\ge 0$ times in the positive direction
(fixing such a direction), then we arrive at $(x,m)$ in $M$. If the
 curve goes through $\Sigma$ $m\ge 0$ times in the 
 negative direction, then we arrive at $(x,-m)$ in $M$.
For $n\ge 3$, the fundamental group of 
$D\setminus \Gamma$ is $\Bbb Z$ and
$M$ is the universal cover of $D\setminus \Gamma$.

For $k=2,3,4, \cdots$, we introduce an equivalence relation
``$\sim_k$'' on $M$ as follows:
$(x,m)$ and $(y,l)$ in $M$ are ``$\sim_k$'' equivalent if
$x=y$ and $m-l$ is an integer multiple
of $k$.  We let
$$
M_k:= M/ \sim_k,
$$
denote the $k-$sheet cover of
$D\setminus \Gamma$, and let
$$
\partial ' M_k
:= \cup_{ i=1}^k (\partial D\times \{i\}).
$$

For $\varphi_1, \cdots, \varphi_k\in C^0(\partial D)$,
it is easy to prove, by Perron$'$s method, that there exists 
$h\in C^\infty(M_k)\cap L^\infty(M_k)\cap C^0(M_k
\cup \partial' M_k)
$ satisfying
\begin{equation}
\left\{
\begin{array}{rll}
\Delta h&=0,&\mbox{on}\ M_k\\
h&=\varphi_i,& \mbox{on}\
\partial D\times \{i\}, \ 1\le i\le k.
\end{array}
\right.
\label{3-1}
\end{equation} 

Since $\Gamma$ has zero capacity,
the maximum principle holds
 on $M_k$:  Let $u, v\in L^\infty(M_k)$ satisfy
$\Delta u\ge 0\ge \Delta v$  in $M_k$
and $\liminf_{dist(y, \partial' M_k)}(u(y)-v(y))\le 0$, then
$u\le v$ in $M_k$.

Let $\bar h\in C^\infty(D)\cap C^0(\overline D)$ be the solution to
$$
\left\{
\begin{array}{rll}
\Delta \bar h&=0, & \mbox{in}\ D,\\
\bar h&=\frac 1k \sum_{ i=1}^k \varphi_i,&
\mbox{on}\ \partial D.
\end{array}
\right.
$$

It was proved by Caffarelli \cite{C5}, under some mild additional 
regularity assumption on $\Gamma$ (e.g. $\Gamma$ is $C^1$), that 
$$
\lim_{ x\to \bar x}h(x, m)
=\bar h(\bar x),\quad
\forall \bar x\in \Gamma, \ 1\le m\le k,
$$
and, for some $0<\alpha<1$ and $C>0$, that 
$$
|h(x,m)-\bar h(\bar x)|\le C|x-\bar x|^\alpha,
\quad \forall\ (x,m)\in M_k.
$$

It follows, by the maximum principle, that
(\ref{3-1}) has a unique bounded
solution $h$.

\medskip

Let $\varphi_1, \cdots, \varphi_k\in C^0(\partial D)$, and let
$f\in C^0(M_k)$ satisfy, for some positive constants $a$ and $b$,
\begin{equation}
a\le f\le b\qquad \mbox{on}\ M_k.
\label{5-1}
\end{equation}

We consider the following Monge-Amp\`ere
equation on $M_k$ with Dirichlet boundary condition:
\begin{equation}
\left\{
\begin{array}{rll}
\det(D^2u)&=f,&\mbox{on}\ M_k,\\
u&=\varphi_i,&
\mbox{on}\
\partial D\times \{i\}, \
1\le i\le k.
\end{array}
\right.
\label{2-1}
\end{equation}

\begin{thm}  Let $M_k$ be as above, $k=2,3,4,\cdots$,
$\varphi_1, \cdots, \varphi_k\in C^0(\partial D)$,
$h$ be the bounded solution of (\ref{3-1}), and let $f\in C^0(M_k)$ satisfy
(\ref{5-1}) for some positive constants $a$ and $b$.
Then (\ref{2-1}) has at least one bounded locally convex viscosity solution
$u$ satisfying $u\le h$ on $M_k$.
\label{thm1}
\end{thm}

\noindent{\bf Proof.}\ Let $P(x)$ be   a convex
quadratic polynomial  satisfying 
$$
\det(D^2P)\ge b\quad \mbox{on}\ D,
$$
$$
P< \inf_{M_k}h\qquad
\mbox{on}\  \overline D, 
$$
and let $D'$ be an open set in $D$ containing
$\Sigma$ and satisfying $dist(D', \partial D)>0$.
As in Caffarelli, Nirenberg and Spruck  \cite{CNS}, 
we construct ${\underline u}_i\in C^\infty(D)\cap C^0(\overline D)$,
$1\le i\le k$, which satisfy
$$
\det({\underline u}_i)\ge b, \qquad \mbox{on}\ D,
$$
$$
{\underline u}_i=\varphi_i, \qquad \mbox{on}\ \partial D,
$$
$$
{\underline u}_i<P, \qquad \mbox{on}\  D'.
$$

Define 
$$
{\underline u}(x,m)=\max\{
{\underline u}_m(x), P(x)\},
\qquad x\in D\setminus \Gamma, \ 1\le m\le k.
$$
Then ${\underline u}\in
C^0(M_k \cup \partial' M_k)$ is a locally convex subsolution of
(\ref{2-1})  satisfying
$$
{\underline u}(x,m)=P(x)
\qquad \forall\ x\in D', 1\le m\le k.
$$

Let ${\cal S}$ denote the set of
locally convex functions $v$ in $C^0(M_k\cup \partial'M_k)$
which are  viscosity subsolutions to (\ref{2-1})
satisfying
\begin{equation}
\limsup_{ x\to \bar x}
\max_{ 1\le m\le k}[v(x,m)-h(x,m)]\le 0,
\qquad \forall\
\bar x\in \Gamma.
\label{7-1}
\end{equation}
Clearly ${\underline u}\in {\cal S}$.

Define on $M_k$
$$
u(x,m)=\sup \{v(x,m)\
|\ v\in {\cal S}\},\quad 1\le m\le k.
$$

For every $v\in {\cal S}$,
$$
\Delta v\ge 0\qquad \mbox{on}\ M_k.
$$
By the maximum principle, using (\ref{7-1}), 
$
v\le h$ on $M_k$.
Thus
$ u\le h$ on $M_k$, and 
$u\in C^0(M_k\cup \partial' M_k)$
is a locally convex viscosity solution of (\ref{2-1}).
Theorem \ref{thm1} is established.

\vskip 5pt
\hfill $\Box$
\vskip 5pt

Let ${\cal S}^*$ denote the set of locally convex functions $v$
in $C^0(M_k\cup \partial' M_k)$
which are viscosity solutions to (\ref{2-1})
satisfying (\ref{7-1}).
 Then
$$
u^*(x,m):=\sup\bigg\{ v(x,m)\
|\ v\in {\cal S}^*\bigg\},
\quad x\in D\setminus\Gamma, 1\le m\le k
$$
is the largest element in 
${\cal S}^*$.  Moreover, by the maximum principle,
$$
u^* \le h\qquad \mbox{in}\
M_k.
$$

It is clear from the proof of Theorem \ref{thm1} that 
(\ref{2-1}) has infinitely many solutions.  In the following
we study the existence of solutions to (\ref{2-1})
with prescribed values on $\Gamma$.

\begin{thm}  Let $M_k$, $\varphi_i$, $f$ and
$h$ be as in Theorem \ref{thm1},
and let ${\underline u}\in C^0(M_k\cup \partial' M_k)$
be a locally convex viscosity subsolution of (\ref{2-1}) satisfying
$$
{\underline u}(\bar x):=\lim_{x\to \bar x}
{\underline u}(x,m)\
\mbox{exists, finite, independent of}\
1\le m\le k, \forall\ \bar x\in \Gamma,
$$
and
$$
{\underline u}(\bar x)\le
\liminf_{x\to\bar x}\min_{1\le m\le k}h(x,m),\qquad
 \forall\ \bar x\in \Gamma.
$$
Then there exists a unique locally convex viscosity
solution $u$ of (\ref{2-1}) satisfying
$$
u(\bar x):=\lim_{x\to \bar x}u(x,m)
={\underline u}(\bar x), \qquad
1\le m\le k, \forall\ \bar x\in \Gamma.
$$
\label{thm2}
\end{thm}

\noindent{\bf Proof.}\ Let ${\cal S}$ denote the set of locally convex
viscosity subsolutions $v$ of (\ref{2-1}) in
$C^0(M_k\cup \partial' M_k)$ satisfying
\begin{equation}
\limsup_{ x\to \bar x}v(x,m)
\le {\underline u}(\bar x),
\quad 1\le m\le k,
\forall\ \bar x\in \Gamma.
\label{12-1}
\end{equation}
Clearly ${\underline u}\in {\cal S}$.
Define on $M_k$,
$$
u(x,m):=\sup \bigg\{
v(x,m)\ |\ v\in {\cal S}\bigg\},
\qquad
1\le m\le k, x\in \overline D\setminus \Gamma.
$$
By the maximum principle, 
$$
{\underline u}\le u\le h\qquad
\mbox{on}\ M_k.
$$
It follows that $u\in C^0(M_k\cup \partial' M_k)$
is a locally convex viscosity solution
of (\ref{2-1})  satisfying
$$
\liminf_{x\to \bar x} u(x,m)\ge {\underline u}(\bar x),
\qquad \forall\ 1\le m\le k, \forall\ \bar x\in \Gamma.
$$
Let $v\in {\cal S}$,  $\bar x\in \Gamma$ and
$x\in D\setminus \Gamma$.
Since the Hausdorff measure $H^{n-1}(\Gamma)=0$,
there exist $y_i\to \bar x$ and $\nu_i\to \frac {x-\bar x}{ \|x-\bar x\|}$
such that
$$
\{y_i+t\nu_i\ |\ t\ge 0\}\cap \Gamma=\emptyset.
$$
By the convexity and the boundedness of $v$ on the lifting
of the segment $\{y_i+t\nu_i\ |\ t\ge 0\}\cap D$, we have,
for some constant $C$ independent of $i$,
$$
v(y_i+|x-\bar x|\nu_i, m)
\le \max_{1\le m'\le k}v(y_i, m')
+C|x-\bar x|, \qquad
\forall\ 1\le m\le k.
$$
Sending $i$ to infinity, we have, by (\ref{12-1}),
$$
v(x,m)\le {\underline u}(\bar x)+C|x-\bar x|,\qquad
\forall\ x\in D\setminus \Gamma, \forall\ 1\le m\le k.
$$
It follows that
$$
\limsup_{ x\to \bar x} u(x,m)\le {\underline u}(\bar x).
$$
Theorem \ref{thm2} is established.

\vskip 5pt
\hfill $\Box$
\vskip 5pt

\noindent $\underline { \mbox{Example 1} }.$\
Let $M_k$ and $f$ be as in Theorem \ref{thm1}, and let
$\varphi_1, \cdots, \varphi_k\in C^0(\partial D)$ satisfy
$$
\varphi_i \ge \varphi_1\qquad
\mbox{on}\ \partial D, 1\le i\le k.
$$
Let $G\in C^0(\overline D)$ be a convex function
satisfying, in the viscosity sense,
$$
\left\{
\begin{array}{rll}
\det(D^2 G)&\ge f&\mbox{in}\ D,\\
G&=\varphi_1&\mbox{on}\ \partial D.
\end{array}
\right.
$$
Then (\ref{2-1}) has a unique locally convex viscosity
solution 
 with $G$ as the prescribed value on
$\Gamma$.

\vskip 5pt

Let $D'$ be a convex open set containing
$\Gamma$ satisfying $\overline {D'}\subset D$.  We follow
\cite{CNS}
to construct convex ${\underline u}_2', \cdots,
{\underline u}_k'\in C^\infty(D)\cap C^0(\overline D)$
satisfying 
$$
\det(D^2 {\underline u}_i')\ge f\qquad \mbox{in}\ D,
2\le i\le k,
$$
$$
{\underline u}_i'=\varphi_i, \qquad \mbox{on}\ \partial D, 
2\le i\le k,
$$
$$
G> {\underline u}_i'\qquad \mbox{in}\
D', 2\le i\le k.
$$

Let ${\underline u}_1'=G$,
$$
{\underline u}_i(x):=\max\{ {\underline u}_1'(x),
{\underline u}_i'(x)\},
\qquad x\in \overline D,
$$
and
$$
{\underline u}(x,m):= {\underline u}_m(x), \qquad
x\in D\setminus \Gamma, 1\le m\le k.
$$
Then  ${\underline u}\in C^0(M_k\cup \partial' M_k)$
is a locally convex viscosity subsolution
of (\ref{2-1}) satisfying
$$
{\underline u}(x,m)=G(x),
\qquad
\forall \ 1\le m\le k, \ \mbox{for}\ x
\ \mbox{in}\ D\setminus \Gamma\ \mbox{ and close to}\ \Gamma.
$$
As a result, by Theorem \ref{thm2}, we can
solve (\ref{2-1}) with $G$ as the prescribed value on
$\Gamma$.

\vskip 5pt

\section{Classical solutions when
$\Gamma$ is a ``plane curve''}

Solutions given by Theorem \ref{thm1} and Theorem \ref{thm2}
are not necessarily classical solutions.
In this section we study the  existence of classical solutions of
(\ref{2-1}) with value $0$ on $\Gamma$ under
some further hypothesis on $\Gamma$.

Let $\Omega\subset D$ be
two bounded open strictly convex subsets with smooth boundaries,
denoted respectively by $\partial \Omega$ and $\partial D$.
Let $\Sigma$, diffeomorphic to a $(n-1)-$disc,
be the intersection of
$\Omega$ and  a hyperplane in $\Bbb R^n$, and let $\Gamma$ be the boundary
of $\partial \Sigma$.
The fundamental group of
$D\setminus \Gamma$ is
$\pi_1(D\setminus \Gamma)=\Bbb Z$ when $n\ge 3$.  Let
$M=(D\setminus \Gamma)\times \Bbb Z$ and $M_k=M/\sim_k$ be
covering spaces  
of $D\setminus \Gamma$ as in Section 1.
 $\Sigma$ divides $\Omega$ into two open parts, denoted
 as $\Omega^+$ and $\Omega^-$.
Fixing a $x^*\in \Omega^-$, we use the convention that
going through
$\Sigma$  from $\Omega^-$ to $\Omega^+$
denotes the positive direction 
through $\Sigma$.

\begin{thm} Let $M_k$, $k=2,3,4, \cdots$, be as above
and let $f\in C^\infty(M_k)$ satisfy (\ref{5-1}) for some 
positive constants $a$ and $b$.  Then 
there exists some $\beta>0$ such that
for any $\varphi_1, \cdots, \varphi_k\in C^\infty(\partial D)$
satisfying 
\begin{equation}
\varphi_i>\beta\qquad \mbox{on}\ \partial D,
1\le i\le k,
\label{18-1}
\end{equation}
there exists a unique locally convex $u\in 
C^\infty(M_k\cup \partial ' M_k)$ satisfying
(\ref{2-1}) and
\begin{equation}
\lim_{x\to \bar x} u(x,m)=0,
\qquad \forall\ 1\le m\le k, \forall\ \bar x\in \Gamma.
\label{18-2}
\end{equation}
\label{thm3}
\end{thm}

\noindent{\bf Proof of Theorem \ref{thm3}.}\ 
Let $v\in C^\infty(\overline \Omega)$ be 
a convex function satisfying 
$$
\left\{
\begin{array}{rll}
\det(D^2 v)&\ge b\qquad\mbox{in}\ \Omega,
\\
v&=0\qquad \mbox{on}\ \partial \Omega.
\end{array}
\right.
$$
By Lemma \ref{A.1}, there exists $\bar x(\xi)\in \Bbb R^n$ for every $\xi\in
\partial \Omega$ such that
$$
w_\xi(x):=\frac 12\left(|x-\bar x(\xi)|^2-|\xi-\bar x(\xi)|^2\right)
< b^{-\frac 1n}v(x),\quad
\forall\ x\in \overline \Omega\setminus \{\xi\}.
$$
Moreover $\sup_{ \xi\in \partial \Omega}|\bar x(\xi)|<\infty$.

Set
$$
V(x)=\left\{
\begin{array}{ll}
v(x),& x\in \overline \Omega,\\
\sup_{ \xi\in \partial \Omega} \left(b^{\frac 1n}w_\xi(x)\right),&
x\in \overline D\setminus \overline \Omega.
\end{array}
\right.
$$
Then
$V\in C^0(\overline D)$ is a 
convex viscosity subsolution to
$$
\det(D^2 V)=b\qquad \mbox{in}\ D.
$$
Let 
$$
\beta=\max_{\partial D}V.
$$
With this value of  $\beta$, and for any
 $\varphi_1, \cdots, \varphi_k\in C^\infty(\partial D)$
satisfying (\ref{18-1}), we can construct as in \cite{CNS}
 convex ${\underline u}_1',
\cdots, {\underline u}_k'\in C^\infty(D)\cap C^0(\overline D)$
satisfying 
$$
\det(D^2 {\underline u}_i')\ge b\qquad \mbox{in}\ D, 1\le i\le k,
$$
$$
{\underline u}_i'=\varphi_i\qquad \mbox{on}\
\partial D, 1\le i\le k,
$$
$$
V>{\underline u}_i'\qquad \mbox{in}\ \overline \Omega,
1\le i\le k.
$$
Let
$$
{\underline u}_i(x)=\max\{
V(x), {\underline u}_i'(x)\},
\qquad x\in D, 1\le i\le k,
$$
and
$$
{\underline u}(x,m)={\underline u}_m(x),\qquad
x\in D\setminus \Gamma, 1\le m\le k.
$$
Clearly ${\underline u}$ is a locally convex viscosity
subsolution of (\ref{2-1}), ${\underline u}
\in C^0(M_k\cup \partial' M_k)$ and
$$
{\underline u}(x,m)=V(x),\qquad
\forall\ 1\le m\le k, \ \mbox{for}\
x\ \mbox{in}\
D\setminus \Gamma\
\mbox{close to}\ \Gamma.
$$
In particular
$$
{\underline u}(\bar x)=\lim_{ x\ \to \bar x} {\underline u}(x,m)
=V(\bar x)=0,\qquad
\forall\ 1\le m\le k, \forall\ \bar x\in \Gamma.
$$
By Theorem \ref{thm2}, there exists a unique locally convex viscosity solution $u$
of (\ref{2-1}) satisfying (\ref{18-2}).

To complete the proof of Theorem \ref{thm3}, we need to show that 
$u\in C^\infty(M_k\cup \partial' M_k)$.
This follows from the regularity theory
developed by the first author in \cite{C1}, 
\cite{C2} and \cite{C3}.
Indeed if $u$ does not belong to
$C^\infty(M_k\cup \partial' M_k)$, then,
by theorem 1 in \cite{C1} and
theorem 2 in \cite{C2}, there must be a line in $M_k$
on which $u$ is linear.  By the arguments in the proof of
corollary 4 in \cite{C1},
 the line  can not
hit $\partial' M_k$.  So 
this singular line $\gamma(t)$,
$0<t<1$, must be the lifting of
$\{tx^{(1)}+(1-t) x^{(2)}\ |\
0<t<1\}$ for some $x^{(1)}, x^{(2)}\in \Gamma$,
$x^{(1)}\ne x^{(2)}$.
By 
(\ref{18-2}), $u(\gamma(t))=0$ for all
$0<t<1$ and therefore $u=0$
on the $(n-1)$ dimensional disc spanned by
$\{\gamma(t)\}_{ 0<t<1}$ and $\Gamma$.  This 
violates the theorem in \cite{C3}.
Theorem \ref{thm3} is established.

\vskip 5pt
\hfill $\Box$
\vskip 5pt

\section{Infinitely valued solutions with exponentially growing 
right hand side}

This section is motivated by the following $2-$d example:
In the spirit of $arg(z)$ that grows by a constant every time
we go around the origin, we 
construct a solution of the Monge-Amp\`ere equation that grows by a 
factor every time
we go around the origin. In dimension $2$ we get the solution
$$
u(r, \theta)=r^2 e^{\lambda \theta}
$$
that satisfies
$$
\det(D^2 u) =
2(2\lambda^2 -2\lambda+1)e^{\lambda \theta}.
$$

We do now a similar construction in $\Bbb R^n$.
For $D, \Omega, \Gamma, M$  as in the last  section, we
study in this section multi-valued solutions to
Monge-Amp\`ere equations on $M$.

Let $S\in C^\infty(M)$ satisfy
$$
S(x,k)=S(x,k-1)+1,
\qquad \forall\ x\in D\setminus \Gamma, \forall \
k\in \Bbb Z,
$$
and
$$
c:
=\sup_{ |k|\le 2, x\in \overline \Omega} e^{ S(x,k)} < \infty.
$$

We use notation
$$
\partial' M= \cup_{i=-\infty}^{\infty} (\partial D\times \{i\}).
$$

\begin{thm} Let $D, \Omega, \Gamma, M, S$ be as above.
Then there exists some positive constant $\beta$ such that for any
$\varphi\in C^\infty(\partial D)$ satisfying
\begin{equation}
\varphi>\beta \qquad\mbox{on}\
\partial D,
\label{22-2}
\end{equation}
there exists a locally convex $u\in
C^\infty(M\cup \partial' M)$
satisfying
\begin{equation}
\det(D^2u)=e^S,\qquad \mbox{in}\
M,
\label{23-1}
\end{equation}
\begin{equation}
u(x,k)=e^{ \frac 1n} u(x, k-1),
\qquad \forall\ (x,k)\in M,
\label{23-2}
\end{equation}
\begin{equation}
\lim_{x\to \bar x}u(x,m)=0,
\qquad \forall\ \bar x\in \Gamma,
\ \forall\ 1\le m\le k,
\label{23-11}
\end{equation}
\begin{equation}
u(x,k)=e^{ \frac kn } \varphi(x),\qquad
\forall\ (x,k)\in \partial ' M.
\label{23-3}
\end{equation}
\label{thm4}
\end{thm}

\noindent {\bf Proof of Theorem \ref{thm4}.}\ 
Let $\xi\in C^\infty(\overline \Omega)$ be a  convex
function satisfying
$$
\left\{
\begin{array}{rll}
\det(D^2\xi)&\ge c,& \mbox{in}\ \Omega,\\
\xi&=0,& \mbox{on}\
\partial \Omega.
\end{array}
\right.
$$

 As in the proof of Theorem \ref{thm3}, we can extend $\xi$ to 
a convex $\widetilde\xi
\in C^0(\overline D)$ which satisfies
in the viscosity sense
$$
\det(D^2\widetilde \xi)\ge c\qquad \mbox{in}\
D.
$$

Let
$$
\beta:=\max_{\partial D}\widetilde \xi.
$$
With this value of $\beta$, for any $\varphi\in C^\infty(\partial D)$
satisfying (\ref{22-2}), we construct, as in
\cite{CNS}, some convex
$\eta'\in C^\infty(D)\cap C^0(\overline D)$ which satisfies
$$
\det(D^2\eta')\ge c\qquad \mbox{in}\ D,
$$
$$
\eta'=\varphi\qquad \mbox{on}\ \partial D,
$$
$$
\eta' < \widetilde \xi\qquad \mbox{on}\ \overline \Omega.
$$
Let
$$
\eta(x):=\max\{\eta'(x), \widetilde \xi(x)\},
\qquad x\in \overline D.
$$
Then $\eta\in C^0(\overline D)$ is a locally convex function satisfying
\begin{eqnarray*}
\eta &=&\varphi\qquad \mbox{on}\
\partial D,\\
\eta&=& \widetilde \xi
\qquad \mbox{in an open neighborhood
of }\ \overline \Omega,
\end{eqnarray*}
and, in the viscosity sense,
$$
\det(D^2 \eta)\ge c\qquad \mbox{in}\ D.
$$
In particular,
$$
\eta=\xi\qquad\mbox{on}\
\overline \Omega,
$$
$$
\eta=0\qquad \mbox{on}\
\partial \Omega,
$$
$$
\eta<0\qquad \mbox{in}\
\Omega.
$$

Define, for $k\in \Bbb Z$,
$$
{\underline u}(x,k)=
\left\{
\begin{array}{rl}e^{  \frac {k-1}n } \eta(x), &\qquad
x\in \Omega^+,\\
e^{\frac kn}\eta(x),&\qquad x\in
\overline D\setminus \overline {\Omega^+}.
\end{array}
\right.
$$
It is not difficult to see that this
extends to ${\underline u}\in C^0(M\cup \partial' M)$ which
is  locally convex and
satisfies
$$
{\underline u}(x,k)=e^{ \frac 1n} {\underline u}(x,k-1),
\qquad \forall\ (x,k)\in M,
$$
$$
\lim_{x\to\bar x} {\underline u}(x,m)=0,
\qquad \forall\ \bar x\in \Gamma, \ \forall\
1\le m\le k,
$$
$$
{\underline u}(x,k)=e^{ \frac kn}\varphi(x),
\qquad \forall\ (x,k)\in \partial D\times \Bbb Z,
$$
and, in the viscosity sense,
$$
\det(D^2 {\underline u})\ge e^S\qquad \mbox{in}\
M.
$$

Let ${\cal S}$ denote the set of 
locally convex functions $v\in C^0(M\cup \partial'M)$
satisfying 
$$
v(x,k)=e^{ \frac 1n} v(x, k-1), \qquad
\forall\ (x,k)\in M,
$$
$$
\lim_{x\to \bar x}v(x,m)=0, \qquad
\forall\ \bar x\in \Gamma, \ \forall\
1\le m\le k,
$$
$$
v(x,k)= e^{ \frac kn}\varphi(x),\qquad
\forall\ (x,k)\in \partial' M,
$$
and, in the viscosity sense,
$$
\det(D^2 v)\ge e^S\qquad \mbox{in}\ M.
$$

Let $B_1, B_2, B_3, \cdots$ be open balls in
$D\setminus \Gamma$ such that
$$
D\setminus \Gamma =\cap_{ i=1}^\infty\cup _{j=i}^\infty B_j,
$$
i.e. every point in $D\setminus \Gamma$ belongs
to infinitely many balls.

For any $v\in {\cal S}$, and for any 
open ball $B\subset D\setminus \Gamma$, we define
$T_Bv$ as follows:
The lifting of $B$ into $M$ is the union of
infinite disjoint balls,
denoted as
$\{ B^{ (m) } \}_{ m=-\infty }^{\infty}$.
We keep $T_B v$ the same as $v$ outside
$\cup _{ m=-\infty}^{\infty}B^{(m)}$, while in each
$B^{(m)}$ we replace $v$ by the solution
of
\begin{equation}
\left\{
\begin{array}{rll}
\det(  D^2 (T_Bv))&=& e^{ S},
\qquad \mbox{in}\ B^{(m)},\\
(T_Bv)&=& v,\qquad \mbox{on}\  \partial B^{(m)}.
\end{array}
\right.
\label{ball}
\end{equation}
It is not difficult to see that $T_Bv\in
{\cal S}$, and $T_Bv\ge v$ in $M$.  Let
$B_1', B_2', B_3', B_4', \cdots$ be a sequence 
of balls defined by $B_1'=B_1$,
$B_2'=B_2$, $B_3'=B_1$, $B_4'=B_2$, $B_5'=
B_3$, $B_6'=B_1$, $B_7'=B_2$,
$B_8'=B_3$, $B_9'=B_4$, $B_{10}'=B_1, \cdots$, and
let $v_0= {\underline u}$ and
$v_i=T_{B_i'}v_{i-1}$ for $i=1,2,3,\cdots$.
Thus we have defined a sequence of functions
$\{v_i\}$ in ${\cal S}$ which satisfy
$$
v_0\le v_1\le v_2\le v_3\le \cdots\qquad \mbox{in}\ M.
$$

For $x\in \partial D$, let $\nu(x)$ denote
the unit inner normal of $\partial D$
at $x$.  We will show below that there exist
some positive constants $\epsilon$ and $C$
such that
\begin{eqnarray}
v_i(x+t\nu(x), k)&\le & v_i(x,k)+Ce^{ \frac kn}t\nonumber\\
&=& e^{ \frac kn}+C e^{ \frac {k+2}n } t,
\qquad \forall\ x\in \partial D, k\in \Bbb Z, 0<t<\epsilon,
\label{30-2}
\end{eqnarray}
and
\begin{equation}
v_i(x,k)\le C e^{ \frac kn},
\qquad \forall\ x\in D\setminus \Gamma, \forall\ k\in \Bbb Z.
\label{30-3}
\end{equation}

Since $\Gamma$ is closed, there exists $\epsilon>0$ such that 
$dist(\Gamma, \partial D)>\epsilon$.  For any
$x\in \partial D$, since $H^{n-1}(\Gamma)=0$, there exists
$|\nu_l|=1$, $\nu_l\to \nu(x)$, such that
$$
\{x+t \nu_l\ |\
t>0\}\cap \Gamma=\emptyset.
$$
Let $t_l>\epsilon$, 
$x+t_l \nu_l\in \partial D$ and let 
 $\gamma(t), 0\le t\le t_l$, $\gamma(0)=(x,k)$, be the lifting
 of $\{x+t\nu_l\ |\ 0\le t\le t_l\}$ to $M$,
then $v_i(\gamma(t))$ is a convex 
function
for $t\in [0, t_l]$.  Since $\Gamma\in \partial \Omega$ and
 $\Omega$ is strictly convex, the segment $\{x+t\nu_l\ |\ 0\le t\le t_l\}$
  can intersect  $\Sigma$ at most once.
 Therefore,
for $|k'-k|\le 1$,
\begin{eqnarray*}
v_i(x+t\nu_l, k)&\le &
v_i(x,k)+\left(\frac{ v_i(x+t_l\nu_l, k')-v_i(x,k)}{t_l}\right)t\\
&=& v_i(x,k)+\left( \frac{ e^{ \frac {k'}n }\varphi(x+t\nu_l)
-e^{ \frac kn}\varphi(x)}{t_l}\right) t
\le v_i(x,k)+Ce^{ \frac kn}t.
\end{eqnarray*}
Estimate (\ref{30-2}) is established.

For $x\in D$, $dist(x,\partial D)<\epsilon$,
and $k\in \Bbb Z$, we deduce from (\ref{30-2})
that
$$
v_i(x,k)\le C e^{ \frac kn}.
$$
Since $H^{n-1}(\Gamma)=0$, for any $x\in D\setminus \Gamma$ with
$dist(x, \partial D)>\epsilon$, there exists
$|\nu|=1$ such that
$$
\{x+t\nu \ |\
t\in \Bbb R\}\cap \Gamma =\emptyset.
$$
Let $t^-<0<t^+$ satisfy
$x+t^{\pm}\nu\in \partial D$.
Let $\gamma(t), t^-\le t\le t^+$, $\gamma(0)=(x,k)$, be the
lifting of $
\{x+t\nu \ |\
t^-\le t\le t^+\}$ into $M$. As before,
$v_i(\gamma(t^{\pm}))$ are bounded from above by $Ce^{ \frac kn}$.  
Thus by the convexity of $v_i(\gamma(t))$ in $t$, $v_i(x)$ is bounded from 
above by $C e^{ \frac kn}$.  Estimate
(\ref{30-3}) is established.  
With (\ref{30-2}) and (\ref{30-3}), 
and some standard arguments, $v_i$ monotonically
converge to some locally convex 
$u\in C^0(M\cup \partial' M)$ which satisfy (\ref{23-2}), (\ref{23-3}), 
and, in the viscosity
sense, (\ref{23-1}).  Using some arguments similar to those in the proof of
Theorem \ref{thm2}, we see that $u$ satisfies (\ref{23-11}). 
The smoothness of $u$ follows from  
the regularity theory of the first
author  as  used in the proof of Theorem \ref{thm3}.
Theorem \ref{thm4} is established.

\vskip 5pt
\hfill $\Box$
\vskip 5pt

\section{Global finitely valued solutions}
 We present here existence results
 closely related to theorem 1.7 in \cite{CL}.
Let $\Omega$, $\Gamma$, $\Sigma$ be as at the beginning of  Section 2,
and we take $D$ to be $\Bbb R^n$ instead of a bounded strictly convex
open set.  W restrict to $n\ge 3$.
For $k=2,3, \cdots$, we define $M$ and
$M_k$ as at the beginning of Section 2
with $D$ replaced by $\Bbb R^n$.

Let
$$
{\cal A}=\{A\ |\
A\ \mbox{is real}\ n\times n\ \mbox{symmetric
positive definite matrix with}\
\det(A)=1\}.
$$

For $f\in C^0(M_k)$ satisfying, for some positive constants
$a$ and $b$,
\begin{equation}
a\le \inf_{ M_k }f\le \sup _{M_k}f\le b,
\label{f1}
\end{equation}
and
\begin{equation}
\{ f\ne 1\}\ \mbox{is compact},
\label{f2}
\end{equation}
we consider
\begin{equation}
\det(D^2 u)=f\qquad \mbox{on}\ M_k.
\label{f3}
\end{equation}

\begin{thm} For $n\ge 3$, $k\ge 2$, let $M_k$, $\Gamma$ be as above,
and let $f\in C^0(M_k)$ satisfy (\ref{f1}) and (\ref{f2}) for
some positive constants $a$ and $b$.
Then for any  $c_m\in \Bbb R$, $b_m\in \Bbb R^n$
and $A_m\in {\cal A}$, $1\le m\le k$, there exists
some $\beta_*\in \Bbb R$ such that for any $\beta>\beta_*$
there exists a unique locally convex viscosity solution  $u\in 
C^0(M_k)$ of (\ref{f3}) which satisfy
\begin{equation}
\limsup_{ |x|\to \infty} \left( |x|^{n-2}
\bigg| u(x,m)- [\frac 12 x' A_m x+b_m \cdot x+c_m]\bigg|\right)
<\infty, \qquad \forall\ 1\le m\le k,
\label{u1}
\end{equation}
\begin{equation}
\lim_{ x\to \bar x} u(x,m)=-\beta, \qquad
\forall\ \bar x\in \Gamma, \ \forall\ 1\le m\le k.
\label{u2}
\end{equation}
Moreover $u\in C^\infty(M_k)$ provided that $f\in C^\infty(M_k)$.
\label{thm8}
\end{thm}

\noindent{\bf Proof of Theorem \ref{thm8}.}\   For simplicity  we assume
that $f\equiv 1$.  The general case can be obtained by incorporating
some arguments in \cite{CL}.
Let $\Phi\in C^\infty(\overline \Omega)$ be a convex
function satisfying
$$
\left\{
\begin{array}{ll}
\det(D^2\Phi)>1& \mbox{on}\ \overline \Omega,\\
\Phi=0,& \mbox{on}\ \partial\Omega.
\end{array}
\right.
$$
By Lemma \ref{A.1},
 there exists $\bar x(\xi)\in
  \Bbb R^n$ for every $\xi\in \partial \Omega$ such that
   $$
    w_\xi(x):=
     \frac 12 \bigg(
      |x-\bar x(\xi)|^2-
       |\xi-\bar x(\xi)|^2\bigg)< \Phi(x), \quad \forall\
        x\in \overline \Omega\setminus\{\xi\}.
	 $$
	  Moreover $\sup_{ \xi\in \partial \Omega} |\bar x(\xi)|<\infty$.

	  Define
	  $$
	  V(x)=
	  \left\{
	  \begin{array}{ll}
	  \Phi(x),& x\in \overline \Omega,\\
	  \sup_{\xi\in \partial \Omega} w_\xi(x),& x\in \Bbb R^n\setminus
	  \overline \Omega.
	  \end{array}
	  \right.
	  $$
	  Then $V$ is a convex function satisfying,
	  in the viscosity sense,
	  $$
	  \det(D^2V)\ge 1\qquad\mbox{in}\ \Bbb R^n.
	  $$

Fix some $R_1>0$ such that
$$
\Omega\subset B_{R_1}.
$$
Write
\begin{equation}
\frac 12 x' A_m x+b_m \cdot x+c_m
=\frac 12 |(A_m)^{ \frac 12}x+
(A_m)^{ -\frac 12}b|^2 +c-\frac 12
|(A_m)^{ -\frac 12}b|^2.
\label{ff1}
\end{equation}
Let
$$
R_2:=2 \max_{1\le m\le k} \max_{|x|\le R_1}
|(A_m)^{ \frac 12}x+
(A_m)^{ -\frac 12}b|.
$$
Define, for $a>1$,
$$
w_{m,a}(x):=\inf_{B_{R_2}}V+
\int_{ 2R_2}^{  |(A_m)^{ \frac 12}x+
(A_m)^{ -\frac 12}b| } (s^n+a)^{\frac 1n}ds,
\qquad
0<|x|<\infty.
$$

$w_{m,a}$ satisfies
$$
\det(D^2 w_{m,a}(x))=1\qquad
\forall\ 0<|x|<\infty.
$$

By the definition of $R_2$,
\begin{eqnarray}
w_{m,a}(x)&\le & \inf_{ B_{R_2} }V+ \int_{ 2R_2}^{ R_2/2}
(s^n+a)^{ \frac 1n}ds\nonumber \\
&<& \inf_{  B_{R_2}  }V\le V(x),
\qquad \forall\ 1\le m\le  k,\
\forall\ |x|\le R_1.
\nonumber
\end{eqnarray}

Fixing some $R_3> 3R_2$ satisfying
$$
\min_{ 1\le m\le k} \min_{ |x|=R_3}
|(A_m)^{ \frac 12}x+
(A_m)^{ -\frac 12}b|> 3R_2,
$$
we choose $a_1>1$ such that
$$
w_{m,a}(x)>
 \inf_{ B_{R_2} }V
 +\int_{ 2R_2}^{3R_2} (s^n+a)^{\frac 1n}ds
 >V(x),\   \forall\ |x|=R_3,  \forall \ 1\le m\le k,
\forall\ a\ge a_1.
$$

It is easy to see, in view of (\ref{ff1}),  that
$$
w_{m,a}(x)
= \frac 12 x'A_mx+b_m\cdot x+c_m+\mu(m,a)+O(|x|^{2-n})
\qquad \mbox{as}\  |x|\to\infty,
$$
where $\mu(m,a)$,  monotonic and continuous in $a$ for large $a$, tends
to $\infty$ as $a\to \infty$.

Define, for $a\ge a_1$ and $1\le m\le k$,
$$
{\underline u}_{m,a}(x)=
\left\{
\begin{array}{ll}
\max\{V(x), w_{m,a}(x)\}-\mu(m, a),& |x|\le R_3,\\
w_{m,a}(x)-\mu(m, a),& |x|\ge R_3.
\end{array}
\right.
$$

Then, for $1\le m\le k$,
\begin{equation}
{\underline u}_{m,a}(x)=
 \frac 12 x'A_mx+b_m\cdot x+c_m+O(|x|^{2-n})
 \qquad \mbox{as}\  |x|\to\infty,
 \label{18-11new}
 \end{equation}
 $$
 {\underline u}_{m,a}=-\mu(m, a)\qquad\mbox{on}\ \Gamma,
 $$
 $$
 {\underline u}_{m,a}= V\ \ \mbox{in some open  neighborhood of}\
 \overline \Sigma,
 $$
 and $ {\underline u}_{m,a}$
 is a convex function satisfying, in the viscosity sense,
 $$
 \det(D^2{\underline u}_{m,a})\ge 1\qquad \mbox{in}\ \Bbb R^n.
$$

It is easy to see that there exist continuous functions
$a^{(m)}(a)$, $2\le m\le k$, satisfying
$$
\lim_{a\to\infty} a^{(m)}(a)= \infty
$$
and, for $2\le m\le k$,
$$
\mu(m, a^{ (m) }(a))=\mu(1,a)\ \ \ \mbox{for large }\ a.
$$

Define, with the convention $a^{(1)}(a)=a$, 
$$
{\underline u}_a(x,m)={\underline u}_{m,a^{(m)}(a)}(x),\qquad
\forall\ (x,m)\in M_k.
$$
Then ${\underline u}_a$ is a locally convex function
on $M_k$ satisfying
$$
{\underline u}_a(x,m)=
 \frac 12 x'A_mx+b_m\cdot x+c_m+O(|x|^{2-n})
  \qquad \mbox{as}\  |x|\to\infty,
  $$
  $$
  \lim_{ x\to \bar x} {\underline u}_a(x,m)=
  \mu(1,a),
 \qquad \forall\ \bar x\in \Gamma, \
 \forall\ 1\le m\le k,
 $$
 and, in the viscosity sense,
 $$
 \det(D^2 {\underline u}_a)\ge 1\qquad \mbox{in}\ M_k.
 $$

Next we produce appropriate supersolutions.
Let $R_4$ be defined by
$$
\max_{  1\le m\le k} \max_{ |x|=2R_3 }
|(A_m)^{ \frac 12}x+
(A_m)^{ -\frac 12}b| =R_4,
$$
and let
$$
w_m^+(x)=
\left\{
\begin{array}{ll} 
\displaystyle{
\int_{ R_4}^{|(A_m)^{ \frac 12}x+
(A_m)^{ -\frac 12}b|}
[ s^n -(R_4)^n]^{ \frac 1n}ds,
}
&
|
(A_m)^{\frac 12}x+(A_m)^{-\frac 12}b
|\ge R_4,\\
0,& |
(A_m)^{\frac 12}x+(A_m)^{-\frac 12}b
|<R_4.
\end{array}
\right.
$$
Then $w_m^+\in C^1(\Bbb R^n)\cap C^\infty( \Bbb R^n\setminus B_{R_4})$ is
a convex function satisfying 
\begin{equation}
\det(D^2 w_m^+)(x)=1\qquad
\mbox{for}\
|
(A_m)^{\frac 12}x+(A_m)^{-\frac 12}b
|>R_4, 
\label{A-2-0new}
\end{equation}
\begin{equation}
\nabla  w_m^+(x)=0\qquad \mbox{for}\
|
(A_m)^{\frac 12}x+(A_m)^{-\frac 12}b
|=R_4,
\label{A-2-1new}
\end{equation}
$$
w_m^+(x)=0\qquad \forall\ |x|<\frac 32 R_3,
$$
and, for some $\bar \beta(m)\in \Bbb R$,
$$
 w_m^+(x)= 
 \frac 12 x'A_mx+b_m\cdot x+c_m +\bar \beta(m)
 +O(|x|^{2-n}),
 \qquad \mbox{as}\ |x|\to \infty.
 $$

 Define
 $$
 w^+(x,m):=
 w_m^+(x) -\bar \beta(m), \qquad \forall\ (x,m)\in M_k.
 $$
 Clearly $
  w^+$  satisfies, in the viscosity sense
  $$
  \det(D^2 w^+))\le 1, \qquad
  \mbox{on}\ M_k.
  $$

For $\lambda$ large,
$w^+ + \lambda > {\underline u}_a$ on $M_k$.  Let
$$
\bar \lambda_a:=\inf
\{\lambda>0\ |\  w^+ +\lambda>
{\underline u}_a\
\mbox{on}\ M_k\}.
$$

Fix some $a_2\ge a_1$ such that
\begin{equation}
-\mu(1,a)<-\max_{1\le m\le k} \bar \beta(m)\qquad \forall\ a\ge a_2.
\label{A-4-1new}
\end{equation}

By (\ref{18-11new}), (\ref{A-2-0new}),
(\ref{A-2-1new}) and (\ref{A-4-1new}), no touching of
$w^++\bar \lambda_a$ and ${\underline u}_a$ can occur (see
arguments on page 575 of \cite{CL}).  Thus
$\bar \lambda_a=0$ and
$w^+>{\underline u}_a$ on $M_k$ for all
$a\ge a_2$.

Let ${\cal S}_a$ denote the set of locally convex
functions $v$ on $M_k$ satisfying
$$
v\le  w^+\qquad \mbox{on}\ M_k,
$$
$$
\det(D^2 v)\ge 1\qquad \mbox{on}\ M_k,
$$
$$
\limsup_{ x\to\bar x} v(x,m)
\le -\mu(1,a), \qquad
\forall\ \bar x\in \Gamma, \ \forall\ 
1\le m\le k.
$$
Clearly ${\underline u}_a\in {\cal S}_a$.  Define
$$
u_a(x,m)=\sup\{
v(x,m)\ |\ v\in {\cal S}\},
\qquad \forall\ (x,m)\in M_k.
$$
Using some arguments similar to those in the proof of 
Theorem \ref{thm2}, together with some standard 
arguments, we see that $u_a$,
for $a\ge a_2$, is a locally
convex solution to 
(\ref{f3}) with $f\equiv 1$
satisfying (\ref{u1}) 
and (\ref{u2}) with
$\beta =\mu(1, a)$.
To complete the proof of Theorem \ref{thm8},
we only need to prove that $u\in C^\infty(M_k)$.
This follows from the regularity theory
of the first 
author as used in the proof of Theorem \ref{thm3}.
Indeed the only additional observation is that,
because of (\ref{u1}), there can not be a 
ray to infinity on which $u_a$ is linear. 
 Theorem \ref{thm8} is established.

\section{Infinitely valued solutions
with a triple point}

In this section, only in $\Bbb R^3$, we construct more complex 
multi-valued solutions.
Here, the curve defining
the multiple leaved space is like
a ``Mercedes Benz star''
and each time we cross one of the three
holes, we go into a different copy of
$\Bbb R^3\setminus \Gamma$.  In particular, the origin
is a triple point
where the ``three cuts'' coexist.
We point out that this construction is possible due to
the particular geometry of the Pogorelov singular solution.

 Let  $\Bbb R^3=\{ (x_1, x_2, x_3)\
|\ x_i\in \Bbb R\}$, and let
 $e_1,  e_2, e_3$ be distinct
unit vectors lying in the
$(x_1, x_2)-$plane.  We assume that
\begin{equation}
e_1\cdot e_2>-1,
\quad e_2\cdot e_3>-1,\quad
e_3\cdot e_1>-1.
\label{34-1}
\end{equation}

Let
$$
e_1'=\frac {e_1+e_3}{ 1+ e_1\cdot e_3}, \quad
e_2'=\frac {e_2+e_1}{ 1+ e_2\cdot e_1}, \quad
e_3'=\frac {e_3+e_2}{ 1+ e_3\cdot e_2}, 
$$
and
$$
\ell_1'(x):= e_i' \cdot x, \qquad i=1,2,3.
$$
Clearly
\begin{equation}
\left\{
\begin{array}{rl}
\ell_1'(e_1)&=\ell_1'(e_3)=1,\\
\ell_2'(e_2)&=\ell_2'(e_1)=1,\\
\ell_3'(e_3)&=\ell_3'(e_2)=1.
\end{array}
\right.
\label{36-1}
\end{equation}

Let $D$ be a strictly convex bounded open set containing the
origin with diameter $diam(D)=2$.
Recall the singular solution to Monge-Amp\`ere equation of 
Pogorelov in $3-$dimension:
\begin{equation}
P(x_1, x_2, x_3)=f(x_1)|(x_2, x_3)|^{ \frac 43}
\label{p}
\end{equation}
where $f$ is positive and smooth in, say,
$(-2, 2)$, and
blows up at $x_1=\pm 2$. See, e.g., \cite{C3} for the ODE satisfied by $f$. 
The function $P$ satisfies in the viscosity sense
$$
\det(D^2P)=1.
$$
Moreover $P(x_1, 0,0)=0$ for all $|x_1|<2$.

Let $P_{e_j}(x)$ denotes the Pogorelov solution which vanishes along 
the $e_j-$line,
and let 
$$
h^0(x)=\max_{j}
\{ x\cdot e_j +P_{e_j}(x)\},
\qquad x\in \overline D.
$$
Let 
$$
\beta':=\sup_D\{h^0, \ell_1', \ell_2', \ell_3'\}>0.
$$
For any $\varphi' \in C^0(\partial D)$
satisfying $\min_{ \partial D}\varphi'>\beta'$,
we can construct as before convex $h'\in C^\infty(D)\cap
C^0(\overline D)$ satisfying
$$
\det(D^2 h')>1\qquad \mbox{in}\ D,
$$
$$
h'=\varphi'\qquad \mbox{on}\ \partial D,
$$
$$
h'<0\qquad \mbox{in}\ B_1.
$$

Define
$$
h(x)=\max\{ h^0(x), h'(x)\},\qquad
x\in \overline D.
$$
By (\ref{36-1}) and the fact that $P_{e_j}=0$
along the $e_j-$line, we have,
$\forall\ 0<s<1$,
$$
\ell_1'(se_1)=h(se_1)=s, \quad \ell_1'(se_3)=h(se_3)=s,
$$
$$ 
\ell_2'(se_2)=h(se_2)=s, \quad \ell_2'(se_1)=h(se_1)=s,
$$
$$
\ell_3'(se_3)=h(se_3)=s, \quad \ell_3'(se_2)=h(se_2)=s.
$$
We now consider the convex domain
$$
C_j:=\{x\in D\ |\ h(x)<\ell_j'(x)\}
$$
and let

$$
\Sigma_j=C_j\cap \{x\ |\ x_3=0\},
$$
$$
\Gamma_j:=\partial C_j\cap \{x\ |\ x_3=0\},
$$
$$
\Gamma=\Gamma_1\cup \Gamma_2\cup \Gamma_3,
$$
$$
\Omega_j^+=C_j\cap \{x\ |\ x_3>0\},
\qquad
\Omega_j^-=C_j\cap \{x\ |\ x_3<0\}.
$$

Let $M$ denote the universal cover of $D\setminus \Gamma$.  The fundamental 
group $G$ of $M$ is the free products of three cyclic groups
$G_1$, $G_2$ and $G_3$.  We use $g_i$ to denote
the generator of $G_i$.   Fixing a base point $x^*$ in $D\setminus \Gamma$,
we parameterize points of $M=(D\setminus \Gamma)\times G$ as usual:
For $x\in D\setminus \Gamma$ and  $g= g_1^{i_1}g_2^{i_2}g_3^{i_3}\cdots
g_1^{i_{3l+1}}g_2^{i_{3l+2}}g_3^{i_{3l+3}}\in G$, we use $(x,g)$ to denote the
point of $M$ obtained by a path
 starting from $x^*$,
ending at $x$, and crossing
$\Sigma_1$ $i_1$ times ($i_1=0$ means no crossing,
$i_1>0$ means crossing in the positive direction, i.e., from
$\Omega_1^-$ to $\Omega_1^+$,
$i_1<0$ means crossing in the negative direction), crossing
$\Sigma_2$ $i_2$ times, crossing $\Sigma_3$ $i_3$ times,
crossing $\Sigma_1$ $i_4$ times, crossing $\Sigma_2$
$i_5$ times, $\cdots$, crossing $\Sigma_3$ $i_{3l+3}$ times.
We use notation $\partial' M=\{(x, g)\ |\ x\in \partial D, g\in G\}$.

Let $a_1, a_2, a_3\in \Bbb R$ and let $S$ be a smooth
function defined on $M$ satisfying
\begin{eqnarray}
S(x, g)
&=& S(x, \bar g)+
(i_1+i_4+\cdots+i_{3l+1})a_1
\nonumber\\
&&+
(i_2+i_5+\cdots+i_{3l+2})a_2+
(i_3+i_6+\cdots+i_{3l+3})a_3,
\nonumber
\end{eqnarray}
where $\bar g$ denotes the identity
element of $G$ and $g= g_1^{i_1}g_2^{i_2}g_3^{i_3}\cdots
g_1^{i_{3l+1}}g_2^{i_{3l+2}}g_3^{i_{3l+3}}\in G$.

We will produce in the rest of this section
locally convex viscosity solutions $u$  to
\begin{equation}
\det(D^2 u)=e^S\qquad \mbox{in}\ M,
\label{S2}
\end{equation}
satisfying
\begin{equation}
D^2\bigg(u(x,g)-\gamma(g)u(x, \bar g)\bigg)=0,\qquad \forall\ (x,g)\in M,
\label{S3}
\end{equation}
where 
$$
\gamma(g)=\gamma_1(g)\gamma_2(g)\gamma_3(g),
$$
$$
\gamma_1(g)=
e^{  \frac {a_1}n
(i_1+i_4+\cdots+i_{3l+1})  }, \ \ \gamma_2(g)=
e^{  \frac {a_2}n
(i_2+i_5+\cdots+i_{3l+2}) }, \ \
\gamma_3(g)=
e^{  \frac {a_3}n
(i_3+i_6+\cdots+i_{3l+3})}.
$$

Let $b=e^{  10(|a_1|+|a_2|+|a_3|)}$,
$h^*=bh$, and  $\beta =b\beta'$.
We consider  $\varphi\in C^0(\partial D)$
satisfying 
\begin{equation}
\min_{\partial D}\varphi>\beta.
\label{S4}
\end{equation}
We will first construct ${\underline u}$ on $M$ satisfying
\begin{equation}
\det(D^2 {\underline u})\ge e^S\qquad\mbox{on}\ M,
\label{aa1}
\end{equation}
\begin{equation}
{\underline u}(x,g)=\gamma(g)
 {\underline u}(x,\bar g)
  - \gamma_1(g)
  \ell_1-
  \gamma_2(g)
  \ell_2-
  \gamma_3(g)
  \ell_3, \qquad (x,g)\in M,
  \label{aa2}
  \end{equation}
and
\begin{equation}
{\underline u}(x,g)= \gamma(g) 
\varphi
 - \gamma_1(g)
\ell_1-
\gamma_2(g) \ell_2-
\gamma_3(g)\ell_3, \qquad (x,g)\in \partial' M,
\label{aa3}
\end{equation}
where $\ell_j=b\ell_j'$.

We will use $\partial C_j\cap \{x\ |\ x_3>0\}$ as a cut-off surface,
playing a similar role as $\partial \Omega\cap \partial \Omega^+$ in
the proof of Theorem \ref{thm4}.
We change $h^*-\ell_j$  ($\ell_j=b\ell_j'$) to
$e^{ \frac{a_j}n} (h^*-\ell_j)$ when
crossing $\Sigma_j$ in the positive direction 
into $\Omega_j^+$, so on that leaf, we replace
$h^*$ by $e^{  \frac{a_j}n }(h^*-\ell_j)+\ell_j$.
In general, for $e^{ \frac cn}h^* +\ell$
($\ell$ is some linear function), we change it to
$e^{  \frac{c+a_j}n} (h^*-\ell_j)+e^{\frac cn}\ell_j+\ell$.

Following the above procedure we have defined ${\underline u}$:
$$
{\underline u}(x,\bar g)=
\left\{
\begin{array}{rl}
h^*(x)& x\in \overline D\setminus
(\Omega_1^+\cup \Omega_2^+\cup \Omega_3^+),\\
e^{ \frac {a_j}n }(h^*-\ell_j)+\ell_j&
x\in \Omega_j^+, j=1,2,3,
\end{array}
\right.
$$
and
 ${\underline u}$ satisfies
(\ref{aa1}), (\ref{aa2}) and (\ref{aa3}).

\begin{thm} For $\varphi\in C^0(\partial D)$ satisfying (\ref{S4}),
there exists a unique locally convex viscosity solution
$u\in C^0(M\cup \partial' M)$ to 
(\ref{S2}) satisfying
$$
u(x,g)=\gamma(g)
u(x,\bar g)
  - \gamma_1(g)
    \ell_1-
      \gamma_2(g)
        \ell_2-
	  \gamma_3(g)
	    \ell_3, \qquad (x,g)\in M,
$$
$$
u(x,g)={\underline u}(x,g),\qquad \forall\ (x,g)\in \partial' M,
$$
and
\begin{equation}
\lim_{x\to \bar x}\bigg(u(x,g)- {\underline u}(x,g)\bigg)=0,
\qquad \forall\ \bar x\in \Gamma, g\in G.
\label{S6}
\end{equation}
Consequently $u$ satisfies (\ref{S3}).
\label{thm7}
\end{thm}

\begin{rem}  The above theorem can easily be extended to
$m\ge 3$ unit vectors $e_1, \cdots, e_m$ lying in the
$(x_1, x_2)$-plane satisfying
$$
e_1\cdot e_2>-1, \
e_2\cdot e_3>-1, \cdots, e_{ m-1}\cdot e_m>-1,
\ e_m\cdot e_1>-1.
$$
In fact, 
$\{e_1, \cdots, e_m\}$  do not need to lie exactly in the 
$(x_1, x_2)$-plane.  These can be seen from the proof of Theorem \ref{thm7}
\end{rem}

\noindent {\bf Proof of Theorem \ref{thm7}.}\ 
Let ${\cal S}$ denote the set of locally convex functions
$v\in C^0(M\cup \partial' M)$ satisfying, in the viscosity sense,
$$
\det(D^2v)\ge e^S\qquad \mbox{in}\ M,
$$
$$
\limsup_{x\to \bar x}\bigg(v(x,g)-{\underline u}(x,g)\bigg)
\le 0,\qquad \forall\ \bar x\in \Gamma, \ g\in G,
$$
$$
v(x,g)=\gamma(g)v(x,\bar g)
-\gamma_1(g)\ell_1-\gamma_2(g)\ell_2-\gamma_3(g)\ell_3,
\quad (x, g)\in M,
$$
and
$$
v(x,g)=\gamma(g)\varphi
-\gamma_1(g)\ell_1-\gamma_2(g)\ell_2-\gamma_3(g)\ell_3,
\quad (x, g)\in \partial' M.
$$
Clearly ${\underline u}\in {\cal S}$.
Define
$$
u(x,g)=\sup\ \{v(x,g)\ |\ v\in {\cal S}\},\qquad (x,g)\in M.
$$
Modifying the arguments in the proof of Theorem \ref{thm4},
we see that  $u$ belongs to ${\cal S}$ and satisfies (\ref{S2})
and (\ref{S6}).
The uniqueness of such $u$ follows from standard arguments.
Theorem \ref{thm7} is established.

\vskip 5pt
\hfill $\Box$
\vskip 5pt

\section{Infinitely valued
solutions with constant right hand side}

In this section we construct infinitely valued
solutions with constant right hand side.  The invariance 
here is given by the fact that
$u$ in consecutive leaves differs from
the previous one in an affine transformation.

For $n\ge 3$, we use  
$
\Bbb R^{n-1}=\{x=(x_1, \cdots,  x_{n-1}, 0)\ |\ x_i\in \Bbb R\}
$ to
denote the hyperplane in 
$
\Bbb R^n=\{x=(x_1, \cdots,  x_n)\ |\ x_i\in \Bbb R\}
$.
Let $\Sigma\subset \Bbb R^{n-1}$ be a $(n-1)-$dimensional strictly convex
bounded open set with smooth boundary, and 
 $T$ be a $n\times n$ real matrix satisfying
 $Tx=x$ for all $x\in \Bbb R^{n-1}$, i.e.
$$
T=
\left(
\begin{array}{cccccc}
1&0&0&\cdots&0&\lambda_1\\
0&1&0&\cdots&0&\lambda_2 \\
\cdots&\cdots&\cdots&\cdots&\cdots&\cdots\\
0&0&0&\cdots&1&\lambda_{n-1}\\
0&0&0&\cdots&0&1
\end{array}
\right).
$$
Let $\Gamma$ be the boundary of $\Sigma$ in the $\Bbb R^{n-1}$,
 and let
$$
M=(\Bbb R^n\setminus \Gamma)\times \Bbb Z
$$ be the universal cover of $\Bbb R^n\setminus \Gamma$
with the usual parameterization:
Fixing a point $x^*$ in $\Bbb R^n\setminus \Gamma$,
 and connecting $x^*$ by a smooth curve
 in  $\Bbb R^n\setminus \Gamma$  to a point $x$ in  $\Bbb R^n\setminus \Gamma$.
 If the curve goes through $\Sigma$ $m\ge 0$ times in the positive direction
 (say, increasing $x_3$),  then we arrive at $(x,m)$ in $M$. If the
  curve goes through $\Sigma$ $m\ge 0$ times in the
   negative direction, then we arrive at $(x,-m)$ in $M$.

For $k=2,3,4, \cdots$, we introduce an equivalence 
relation ``$\sim_k$'' on $M$ as follows:
$(x,m)$ and $(y,l)$ in $M$ are ``$\sim_k$'' equivalent if
$x=y$ and $m-l$ is an integer multiple
of $k$.  We let
$$
M_k:= M/ \sim_k,
$$
denote the $k-$sheet cover of
$\Bbb R^n\setminus \Gamma$.

\begin{thm} Let $M$
 and $T$ be as above,  
 $b\in \Bbb R^n$, $c\in \Bbb R$,
 and let $A$ be a real
 symmetric $n\times n$ matrix with $\det(A)=1$.
Then there exists some constant $\beta_*>0$,
which depends only on $T$, $\Gamma$, $A$, $b$ and $c$,
such that
for any $\beta \ge \beta_*$ there exists a unique 
locally convex function
$u\in C^\infty(M)$ satisfying
\begin{equation}
\det(D^2u)=1\qquad \mbox{in}\ M,
\label{T1}
\end{equation}
\begin{equation}
u(x,m)=u(Tx, m-1)\qquad \forall\ (x,m)\in M,
\label{T111}
\end{equation}
\begin{equation}
\lim_{x\to \bar x} u(x,m)=-\beta,\qquad
\forall\ \bar x\in \Gamma, \ \forall\ m\in
{\cal Z},
\label{T4}
\end{equation}
and
\begin{equation}
\limsup_{|x|\to\infty}\bigg( |x|^{n-2}\left| u(x, m)-Q(T^mx)\right|\bigg)
<\infty,
\qquad  \forall\ m\in
{\cal Z},
\label{T3}
\end{equation}
where $Q(x):= \frac 12 x'Ax+b\cdot x+c$.
\label{thm5}
\end{thm}
\begin{rem} It is clear that 
the conclusion of Theorem \ref{thm5} holds with
$M$ replaced by $M_k$, $k=2,3,4, \cdots$.  
\label{rem5-1}
\end{rem}

\noindent{\bf Proof of Theorem \ref{thm5}.}\  By the affine invariance of
the equation, we may assume without loss of generality that
$b=0$, $c=0$ and $A$ is the identity matrix.
Let $\Omega\subset \Bbb R^n$ be a strictly convex bounded open set with
smooth boundary satisfying $\Sigma=\Omega\cap \Bbb R^{n-1}$
and, therefore, $\Gamma=\partial \Omega\cap \Bbb R^{n-1}$.
Let $\Phi\in C^\infty(\overline \Omega)$ be a convex
function satisfying
$$
\left\{
\begin{array}{ll}
\det(D^2\Phi)>1& \mbox{on}\ \overline \Omega,\\
\Phi=0,& \mbox{on}\ \partial\Omega,
\end{array}
\right.
$$
and let
$$
\widetilde \Phi(x):= \Phi(x)+ K|x_n|,\qquad x\in \overline \Omega,
$$
where $K>1$ is some fixed large constant, depending only on $T$ and $\Gamma$,
satisfying, for all $ (x_1, \cdots, x_{n-1}, 0)\in \overline \Sigma$,
that
\begin{equation}
\left\{
\begin{array}{l}
\displaystyle{
\liminf_{s\to 0}
\frac{  \widetilde \Phi(x_1, \cdots, x_{n-1}, s)
-  \widetilde \Phi(x_1, \cdots, x_{n-1}, 0) }{  |s|  } 
>0,
}\\
\displaystyle{
\liminf_{s\to 0}
\frac{  \widetilde \Phi (T(x_1, \cdots, x_{n-1}, s))
-  \widetilde \Phi(T(x_1, \cdots, x_{n-1}, 0)) }{  |s|  }
>0,
}
\end{array}
\right.
\label{K}
\end{equation}
Because of (\ref{K}), $\widetilde \Phi$ also satisfies,
in the viscosity sense,
$$
\det(D^2\widetilde \Phi)>1\qquad \mbox{in}\ \Omega.
$$

By Lemma \ref{A.1},
 there exists $\bar x(\xi)\in
 \Bbb R^n$ for every $\xi\in \partial \Omega$ such that
 $$
 w_\xi(x):= \widetilde \Phi(\xi)+
 \frac 12 \bigg(
 |x-\bar x(\xi)|^2-
 |\xi-\bar x(\xi)|^2\bigg)<\widetilde  \Phi(x), \quad \forall\
 x\in \overline \Omega\setminus\{\xi\}.
 $$
 Moreover $\sup_{ \xi\in \partial \Omega} |\bar x(\xi)|<\infty$.

Define 
$$
V(x)=
\left\{
\begin{array}{ll}
\widetilde \Phi(x),& x\in \overline \Omega,\\
\sup_{\xi\in \partial \Omega} w_\xi(x),& x\in \Bbb R^n\setminus
\overline \Omega.
\end{array}
\right.
$$
Then $V$ is a convex function satisfying,
in the viscosity sense,
$$
\det(D^2V)\ge 1\qquad\mbox{in}\ \Bbb R^n.
$$
Let  $R_1>0$ satisfy 
$$
\Omega \subset  B_{R_1},
$$
and let
$$
R_2= 2\max_{|m|\le 4}\max_{|x|\le R_1} |T^mx|.
$$
We
 consider, for $a>1$,
$$
w_a(x):= \inf_{ B_{ R_2 } }V
+ \int_{2R_2}^{|x|}(s^n+a)^{\frac 1n}ds,
\qquad 0<|x|<\infty.
$$
By the definition of $R_2$,
$$
w_a(T^mx)\le  \inf_{ B_{R_2} }V+ \int_{ 2R_2}^{ R_2/2}
(s^n+a)^{ \frac 1n}ds
< \inf_{  B_{R_2}  }V\le V(x),
\qquad \forall\ |m|\le 4,
\
\forall\ |x|\le R_1.
$$

Fixing some $R_3> 3R_2$ satisfying
$$
\min_{ |m|\le 4}\min_{ |x|=R_3}|T^mx|>3R_2,
$$
we choose some $a_1>1$ such that
$$
w_a(T^mx)
> \inf_{ B_{R_2} }V
+\int_{ 2R_2}^{3R_2} (s^n+a)^{\frac 1n}ds
>V(x),\ \  \forall\ |x|=R_3, \ \forall \ |m|\le 4.
$$

It is easy to see that
$$
w_a(x)=\frac 12 |x|^2+\mu(a)+O(|x|^{2-n})
\qquad\mbox{as}\ |x|\to\infty,
$$
where $\mu(a)$, monotonic and continuous in $a$ for large $a$, tends
to $\infty$ as $a\to \infty$.

Define, for $a\ge a_1$,
$$
{\underline u}_a^{(0)}(x)=
\left\{
\begin{array}{ll}
\max\{V(x), w_a(x)\}-\mu(a),& |x|\le R_3,\\
w_a(x)-\mu(a),& |x|\ge R_3,
\end{array}
\right.
$$
and
$$
{\underline u}_a^{(m)}(x)= {\underline u}_a^{(0)}(T^mx),
\qquad x\in \Bbb R^n\setminus \Sigma, \  m\in {\cal Z}.
$$

Then, for $m\in {\cal Z}$, 
\begin{equation}
{\underline u}_a^{(m)}(x)=\frac 12 |T^mx|^2 +O(|x|^{2-n})
\qquad\mbox{as}\ |x|\to\infty,
\label{18-11}
\end{equation}
$$
{\underline u}_a^{(m)}=-\mu(a)\qquad\mbox{on}\ \Gamma,
$$
$$
{\underline u}_a^{(m)}= V\ \ \mbox{in some open  neighborhood of}\
\overline \Sigma,
$$
and ${\underline u}_a^{(m)}$
is a convex function satisfying, in the viscosity sense,
$$
\det(D^2{\underline u}_a^{(m)})\ge 1\qquad\mbox{in}\ \Bbb R^n.
$$

Define 
$$
{\underline u}_a(x,m)={\underline u}_a^{ (m) }(x),\qquad
\forall\ (x,m)\in M.
$$
Then ${\underline u}_a$ is a locally convex function 
on $M$ satisfying
$$
{\underline u}_a(x,m)=\frac 12 |T^mx|^2
+O(|x|^{2-n}),\quad
\mbox{as}\ |x|\to \infty,
$$
$$
\lim_{ x\to \bar x} {\underline u}_a(x,m)=-\mu(a),
\qquad \forall\ \bar x\in \Gamma, \
\forall\  m\in {\cal Z},
$$
and, in view of (\ref{K}),  in the viscosity sense,
$$
\det(D^2 {\underline u}_a)\ge 1\qquad \mbox{in}\ M.
$$

Let $R_4$ be defined by
$$
\max_{ |m|\le 4}
\max_{ |x|=2R_3 }|T^mx|=R_4,
$$
and let
$$
w(x)=
\left\{
\begin{array}{ll}
\displaystyle{
\int_{ R_4}^{|x|}
[ s^n -(R_4)^n]^{ \frac 1n}ds,
}
&
|x|\ge R_4,\\
0,& |x|<R_4.
\end{array}
\right.
$$
Then $w\in C^1(\Bbb R^n)\cap C^\infty( \Bbb R^n\setminus B_{R_4})$ is
a convex function satisfying
\begin{equation}
\det(D^2 w)=1\qquad
\mbox{on}\ \Bbb R^n\setminus  B_{R_4},
\label{A-2-0}
\end{equation}
\begin{equation}
\nabla w=0\qquad \mbox{on}\
\partial B_{ R_4 },
\label{A-2-1}
\end{equation}
and, for some $\beta \in \Bbb R$,
$$
w(x)=\frac 12 |x|^2 +\beta+O(|x|^{2-n}),
\qquad \mbox{as}\ |x|\to \infty.
$$

Define
$$
\bar w(x, m)= w(T^mx)-\beta, \qquad \forall\ (x,m)\in M.
$$
Clearly $\bar w$ satisfies, in the viscosity sense
$$
\det(D^2 \bar w)\le 1, \qquad 
\mbox{on}\ M.
$$

For $\lambda$ large, 
$\bar w+\lambda>w_a$ on $M$. Let
$$
\bar \lambda_a:=\inf
\{\lambda>0\ |\ \bar w+\lambda>w_a\
\mbox{on}\ M\}.
$$

Fix some $a_2\ge a_1$ such that
\begin{equation}
-\mu(a)<-\beta\qquad \forall\ a\ge a_2.
\label{A-4-1}
\end{equation}

By (\ref{18-11}), (\ref{A-2-0}),
(\ref{A-2-1}) and (\ref{A-4-1}), no touching of
$\bar w+\bar \lambda_a$ and $w_a$ can occur (see
arguments on page 575 of \cite{CL}).  Thus
$\bar \lambda_a=0$ and
$\bar w>w_a$ on $M$ for all
$a\ge a_2$.

Let ${\cal S}_a$ denote the set of locally convex
functions $v$ on $M$ satisfying
$$
v\le \bar w\qquad \mbox{on}\ M,
$$
$$
v(x,m)=v(Tx, m-1)\qquad \forall\ (x,m)\in M,
$$
$$
\det(D^2 v)\ge 1\qquad \mbox{on}\ M,
$$
$$
\limsup_{ x\to\bar x} v(x,m)
\le -\mu(a), \qquad
\forall\ \bar x\in \Gamma, \ \forall\ 
m\in {\cal Z}.
$$
Clearly ${\underline u}_a\in {\cal S}_a$.  Define
$$
u_a(x,m)=\sup\{
v(x,m)\ |\ v\in {\cal S}\},
\qquad \forall\ (x,m)\in M.
$$
Using some arguments similar to those in the proofs of 
Theorem \ref{thm2} and Theorem \ref{thm4}, together with some standard 
arguments, we see that $u_a$,
for $a\ge a_2$, is a locally
convex solution to 
(\ref{T1})
satisfying (\ref{T3}) with
$Q(x)\equiv \frac 12 |x|^2$ and (\ref{T4}) with
$\beta =\mu(a)$.
To complete the proof of Theorem \ref{thm5},
we only need to prove that $u\in C^\infty(M)$.
This follows from the regularity theory
of the first 
author as used in the proof of Theorem \ref{thm3}.
Indeed the only additional observation is that,
because of (\ref{T3}), there can not be a 
ray to infinity on which $u_a$ is linear.  Theorem \ref{thm5} is established.

\section{Infinitely valued solutions with constant
right hand side and a triple point}

In this section we construct, only in $\Bbb R^3$,
infinitely valued solutions with constant
right hand side and a triple point by combining the arguments in
Section 5 and 6.
As in Section 5, the construction relies on the
geometry of the Pogorelov singular solution.

In  $\Bbb R^3=\{ (x_1, x_2, x_3)\
|\ x_i\in \Bbb R\}$,  let
 $e_1,  e_2, e_3$ be distinct
 unit vectors lying in $ \Bbb R^2:=\{(x_1, x_2, 0)\ |\
 x_i\in \Bbb R\}$
  which satisfy (\ref{34-1}).
With the Pogorelov singular solution in (\ref{p}), we
let $P_{e_j}$ denote the Pogorelov solution which vanishes 
along the $e_j-$line, and we
define
$$
h(x)=
\max_{j}
\{ x\cdot e_j +P_{e_j}(x)\},
\qquad |x|\le \frac 32.
$$
We let
$$
\widetilde h(x):=h(x)+K|x_3|,
\qquad |x|\le \frac 32,
$$
and we will fix some large constant $K$ below.

Let $T_1, T_2, T_3$ be $3\times 3$ real matrices satisfying
$T_ix =x$ for all $x\in \Bbb R^2$.
We now fix some large positive constant $K$ such that
for all $|(x_1, x_2, 0)|\le \frac 32$, for
$i=1,2, 3$, and for all
$m=0, \pm 1, \pm 2$, we have
\begin{equation}
\liminf_{s\to 0}
\frac{  \widetilde h(T_i^m(x_1, x_2,s))
-  \widetilde h(T_i^m(x_1, x_2,  0)) }{  |s|  }
>0,
\label{Knew}
\end{equation}

Let $b\in \Bbb R^n$, $A$ be a real symmetric $3\times 3$ matrix,
then there exists
some  $c_*$, which depends on
 $A$, $b$ and 
 $\tilde h$, such that
 for all $c>c_*$
 we can construct, as in
Section 6, a convex function $W$ on $\Bbb R^n$, satisfying,
$$
W(x)=\widetilde h(x),\qquad \forall\ |x|\le \frac 54,
$$
$$
W(x)=\frac 12 x'Ax+b\cdot x+c+O(\frac 1{|x|}),\qquad \mbox{as}\
|x|\to \infty,
$$
and,  in the viscosity sense, 
$$
\det(D^2 W)\ge 1\qquad \mbox{in}\ \Bbb R^n.
$$

Let $\Gamma_1, \Gamma_2, \Gamma_3$ be three bounded convex
curves lying in $\Bbb R^2$ satisfying
$$
\Gamma_1\cap \{ x\in \Bbb R^2\ |\ 
|x|<1\}= \{se_1\ |\ 0\le s<1\}\cup  \{se_2\ |\ 0\le s<1\},
$$
$$
\Gamma_2\cap \{ x\in \Bbb R^2\ |\
|x|<1\}= \{se_2\ |\ 0\le s<1\}\cup  \{se_3\ |\ 0\le s<1\},
$$
and
$$
\Gamma_3\cap \{ x\in \Bbb R^2\ |\
|x|<1\}= \{se_3\ |\ 0\le s<1\}\cup  \{se_1\ |\ 0\le s<1\}.
$$
We set
$$
\Gamma=\Gamma_1\cup \Gamma_2\cup \Gamma_3.
$$

Let $M$ denote the universal cover of
$\Bbb R^3\setminus \Gamma$.  
 The fundamental
 group $G$ of $M$ is the free products of three cyclic groups
 $G_1$, $G_2$ and $G_3$.  We use $g_i$ to denote
 the generator of $G_i$.   Fixing a base point $x^*$ in $\Bbb R^3
 \setminus \Gamma$,
 we parameterize points of $M=(\Bbb R^3\setminus \Gamma)\times G$ as usual:
 For $x\in \Bbb R^3\setminus \Gamma$ and  $g= g_1^{i_1}g_2^{i_2}g_3^{i_3}\cdots
 g_1^{i_{3l+1}}g_2^{i_{3l+2}}g_3^{i_{3l+3}}\in G$, we use $(x,g)$ to denote the
 point of $M$ obtained by a path
  starting from $x^*$,
  ending at $x$, and crossing
  $\Sigma_1$ $i_1$ times ($i_1=0$ means no crossing,
  $i_1>0$ means crossing in the positive direction, i.e., from
  $\Omega_1^-$ to $\Omega_1^+$,
  $i_1<0$ means crossing in the negative direction), crossing
  $\Sigma_2$ $i_2$ times, crossing $\Sigma_3$ $i_3$ times,
  crossing $\Sigma_1$ $i_4$ times, crossing $\Sigma_2$
  $i_5$ times, $\cdots$, crossing $\Sigma_3$ $i_{3l+3}$ times.

For  $g= g_1^{i_1}g_2^{i_2}g_3^{i_3}\cdots
 g_1^{i_{3l+1}}g_2^{i_{3l+2}}g_3^{i_{3l+3}}\in G$, we use
notation
$$
T(g)= T_3^{i_{3l+3}}
T_2^{i_{3l+2}} T_1^{i_{3l+1}} \cdots T_3^{i_3}T_2^{i_2}T_1^{i_1}.
$$
Now we define a function ${\underline u}$ on $M$ by setting,
for all $(x,g)\in M$,
$$
{\underline u}(x,g)=
W\big(T(g)x\big).
$$
It is clear that 
${\underline u}$ is a locally convex function on
$M$ satisfying
$$
{\underline u}(x,\bar g)\equiv 
W(x),
$$
$$
{\underline u}(x, g)\equiv
{\underline u}(T(g)x, \bar g),
$$
and, in the viscosity sense,
$$
\det(D^2 {\underline u})\ge 1\qquad \mbox{in}\ M.
$$
Clearly ${\underline u}$ satisfies
$$
\limsup_{|x|\to \infty}
\bigg(|x||{\underline u}(x)-Q(T(g)x)|\bigg)<\infty,\qquad \forall\ 
g\in G,
$$
where $Q(x):=\frac 12 x'Ax +b\cdot x+c$.

Modifying the construction of the super solution $\bar w$ in
Section 6, and increasing the value of $c_*$ if necessary (recall
that $c>c_*$),
we may construct a locally convex function $\overline u$ on $M$
satisfying
$$
\overline u\ge {\underline u}\qquad \mbox{on}\ M,
$$
$$
\limsup_{|x|\to \infty}
\bigg(|x||{\overline u}(x)-Q(T(g)x)|\bigg)<\infty,\qquad \forall\ 
g\in G,
$$
$$
\overline u(x,g)\equiv \overline u(T(g)x, \bar g)\qquad
\forall\ (x,g)\in M,
$$
and, in the viscosity sense,
$$
\det(D^2 \overline u)\le 1\qquad \mbox{on}\ M.
$$

Let ${\cal S}$ denote the set of locally convex functions
$v$ on $M$ satisfying,
$$
v\le \overline u\qquad \mbox{in}\ M,
$$
$$
v(x,g)=v(T(g)x, \bar g)\qquad \forall\ (x,g)\in M,
$$
$$
\limsup_{x\to \bar x}
\bigg( v(x,g)-{\underline u}(x,g) \bigg)
\le 0\qquad\forall\ \bar x\in \Gamma, \ \forall\ g\in G,
$$
and, in the viscosity sense,
$$
\det(D^2v)\ge 1\qquad \mbox{in}\ M.
$$

Clearly ${\underline u}\in {\cal S}$.  Define
$$
u(x,g)=
\sup\{
v(x,g)\ |\ v\in {\cal S}\},
\qquad \forall\ (x,g)\in M.
$$

\begin{thm}  The above defined $u$ is a locally convex
viscosity solution to
\begin{equation}
\det(D^2u)=1\qquad \mbox{in}\ M,
\label{aaa1}
\end{equation}
satisfying
\begin{equation}
u(x,g)=u(T(g)x, \bar g)\qquad \forall\ (x,g)\in M,
\label{aaa2}
\end{equation}
\begin{equation}
\lim_{x\to \bar x}
\bigg( u(x,g)-{\underline u}(x,g) \bigg)
=0\qquad\forall\ \bar x\in \Gamma, \ \forall\ g\in G,
\label{aaa3}
\end{equation}
and
\begin{equation}
\limsup_{|x|\to \infty}
\bigg(|x||u(x)-Q(T(g)x)|\bigg)<\infty,\qquad \forall\
g\in G.
\label{aaa4}
\end{equation}
\label{thm7-1}
\end{thm}

\noindent{\bf Proof of Theorem \ref{thm7-1}.}\
This theorem follows from some arguments similar to those used
in the proofs of
the theorems in previous sections.
Let us outline the arguments below.
First we let $B_1, B_2, B_3, \cdots$ be defined as 
in the proof of Theorem \ref{thm4}, 
with $D\setminus \Gamma$ replaced by 
$\Bbb R^3\setminus \Gamma$.
For $v\in {\cal S}$, and for any
open ball $B\subset \Bbb R^3\setminus \Gamma$, we define
$T_Bv$ similar to the definition in the proof of
 of Theorem \ref{thm4}, only changing
 $e^S$ in (\ref{ball}) to $1$.
Define $\{v_i\}\subset {\cal S}$
the same as below (\ref{ball}). 
Since 
$v_i\le \overline u$ in $M$, we deduce, using  also the local convexity and the monotonicity property  of
$\{v_i\}$ (recall that  $v_i\le v_{i+1}$ in $M$),   
that  $v_i$ converges in $C^0_{loc}(M)$ to some
locally convex function $u$.  Clearly $u$ satisfies (\ref{aaa2})
and (\ref{aaa4}).
By some standard arguments, $u$ satisfies (\ref{aaa1}) in the viscosity
sense.  
Since
$$
\limsup_{ x\to \bar x}\bigg(v_i(x,g)-{\underline u}(x,g)\bigg)
\le 0\qquad \forall\ \bar x\in \Gamma, \ \forall\ g\in G,
$$ 
we can deduce (\ref{aaa3}) by using arguments similar to those used in 
the proof of Theorem \ref{thm4}.
 Theorem \ref{thm7-1} is established.

\vskip 5pt
\hfill $\Box$                                                                   \vskip 5pt

\section{Appendix}

\begin{lem} Let $D$ be a strictly convex bounded open set in
$\Bbb R^n$, $n\ge 1$, with $C^2$ boundary, and let $\Phi
\in C^2(\overline D)$.  Then there exists
some constant $C$, depending only on
$n, \Phi$ and $D$, such that for every 
$\xi\in \partial D$, there exists $\bar x(\xi)\in 
\Bbb R^n$ satisfying
$$
|\bar x(\xi)|\le C\ \ \mbox{and}\ \ 
w_\xi<\Phi
\
\mbox{on}\ \overline D\setminus \{\xi\},
$$where
$$
w_\xi(x):=
\Phi(\xi)+\frac 12 \bigg(
|x-\bar x(\xi)|^2-
|\xi-\bar x(\xi)|^2\bigg),
\qquad x\in \Bbb R^n.
$$
\label{A.1}
\end{lem}
\noindent{\bf Proof.}\ It follows from
modification of the proof of lemma 5.1 
in \cite{CL}.

\vskip 5pt
\hfill $\Box$                                                                   
\end{document}